\documentclass[10pt]{article}
\usepackage{stmaryrd}
\usepackage{mathrsfs}
\usepackage[centertags]{amsmath}
\usepackage{amsfonts, dsfont}
\usepackage{amssymb}
\usepackage{amsthm}
\usepackage{newlfont}
\usepackage{bezier,  amsxtra, amsbsy, amsgen, amsopn, amstext}

\input xy
\xyoption{all}

\theoremstyle{plain}
\newtheorem{thm}{Theorem}[section]
\newtheorem{cor}[thm]{Corollary}
\newtheorem{lem}[thm]{Lemma}
\newtheorem{prop}[thm]{Proposition}

\theoremstyle{definition}
\newtheorem{defn}[thm]{Definition}
\newtheorem{remark}[thm]{Remark}
\newtheorem*{ack}{Acknowledgments}

\newcommand{\bd}{\begin{defn}}
\newcommand{\ed}{\end{defn}}
\newcommand{\bl}{\begin{lem}}
\newcommand{\el}{\end{lem}}
\newcommand{\bp}{\begin{prop}}
\newcommand{\ep}{\end{prop}}
\newcommand{\bt}{\begin{thm}}
\newcommand{\et}{\end{thm}}
\newcommand{\bc}{\begin{cor}}
\newcommand{\ec}{\end{cor}}
\newcommand{\br}{\begin{remark}}
\newcommand{\er}{\end{remark}}
\newcommand{\bdi}{\begin{diagram}}
\newcommand{\edi}{\end{diagram}}
\newcommand{\beq}{\begin{eqn}}
\newcommand{\eeq}{\end{eqn}}
\newcommand{\ba}{\begin{array}}
\newcommand{\ea}{\end{array}}
\newcommand{\bpf}{\begin{proof}}
\newcommand{\epf}{\end{proof}}

\newcommand{\Q}{\mathds{Q}}
\newcommand{\Zp}{\mathds{Z}_{p}}
\newcommand{\Qp}{\mathds{Q}_{p}}
\newcommand{\al}{\alpha}

\newcommand{\Si}{\Sigma}
\newcommand{\Ga}{\Gamma}

\newcommand{\La}{\Lambda}
\newcommand{\la}{\lambda}

\newcommand{\Op}{\mathcal{O}}

\newcommand{\M}{\mathfrak{M}}

 \DeclareMathOperator{\Sel}{Sel}
\DeclareMathOperator{\Gal}{Gal} \DeclareMathOperator{\Hom}{Hom}
 
 \DeclareMathOperator{\rank}{rank}
\newcommand{\ot}{\otimes}

\newcommand{\ilim}{\displaystyle \mathop{\varinjlim}\limits}
\newcommand{\plim}{\displaystyle \mathop{\varprojlim}\limits}

\newcommand{\coker}{\mathrm{coker}\,}

\newcommand{\cyc}{\mathrm{cyc}}
\newcommand{\cts}{\mathrm{cts}}

\newcommand{\lra}{\longrightarrow}

\newcommand{\ps}[1]{[[ #1 ]]}

\begin{document}

\title{Akashi series, characteristic elements
and congruence of Galois representations}
\author{Meng Fai Lim \footnote{School of Mathematics and
Statistics, Central China Normal University, 152 Luoyu Road,  Wuhan,
Hubei, P.R.China 430079. E-mail: \texttt{limmf@mail.ccnu.edu.cn}}}
\date{}
\maketitle

\begin{abstract} \footnotesize
\noindent  In this paper, we compare the Akashi series of the
Pontryagin dual of the Selmer groups of two Galois representations
over a strongly admissible $p$-adic Lie extension. Namely, we show
that whenever the two Galois representations in question are
congruent to each other, the Akashi series of one is a unit if and
only if the Akashi series of the other is also a unit. We will also
obtain similar results for the Euler characteristics of the Selmer
groups and the characteristic elements attached to the Selmer
groups.



\end{abstract}

\section{Introduction}

Let $p$ be an odd prime. In this paper, we are interested in
studying the Selmer groups attached to a certain class of Galois
representations. This class of Galois representations will contain
those which arises from an abelian variety having good ordinary
reduction at all primes above $p$, and those which arises from
$p$-ordinary modular forms. We will be concerned with comparing the
Selmer groups of two congruent Galois representations. Over the
cyclotomic $\Zp$-extension, such studies were carried out in
\cite{EPW, Gr94, GV, Ha}. One of the motivation behind these studies
lies in the philosophy that the ``Iwasawa main conjecture" should be
preserved by congruences. Naturally, one will like to make an
analogous study over a noncommutative $p$-adic Lie extension and
this has been carried out in \cite{Ch09, CS12, Sh, SS, SS14} to some
extent.

The aim of this paper is to compare the Akashi series (and its
related invariants) of the Selmer groups of the two congruent Galois
representations over a strongly admissible $p$-adic Lie extension
which is not totally real. Our main result is that under appropriate
assumptions, whenever the two Galois representations in question are
congruent to each other, the Akashi series of the Selmer group of
one is a unit if and only if the Akashi series of the other is a
unit (see Theorem \ref{Akashi congruent}). We will also prove
similar statements for the Euler characteristics of the Selmer
groups (see Theorem \ref{Euler congruent}) and for the
characteristic elements attached to the Selmer groups (see Theorem
\ref{char congruent}).

We now give a brief description of the idea behind our proofs. One
of the main ingredient to our proofs is a result of Matsuno
\cite{Mat} which is a fine analysis on the structure of the Selmer
groups of the Galois representations over a cyclotomic
$\Zp$-extension of a non-totally real field (see also \cite{Ch}).
Building on Matsuno's result, we are able to obtain certain fine
results on the structure of the Selmer group over an admissible
$p$-adic Lie extension which is not totally real (see Lemma
\ref{algebra nonzero}). Roughly speaking, the first part of our
results will show that the Selmer group, if nonzero, cannot be ``too
small". Another ingredient we require is an observation on the
structure of the homology of the dual Selmer groups by the author in
\cite{LimMHG}. Here, again roughly speaking, our result will yield
that if the Akashi series of the Selmer group is a unit, then the
Selmer group is ``too small" (see Proposition \ref{algebra prop})
and so by the previous observation, the Selmer group has to be zero.
Building on this observation (see Section \ref{Akashi section}), we
are reduced to showing that the dual Selmer groups have the same
ranks over the Iwasawa algebra of $\Gal(F_{\infty}/F^{\cyc})$. We
will prove this equality of ranks under certain appropriate
assumptions (see Proposition \ref{Akashi congruent lemma}) and hence
establish our main results. It is here that we will require the
second result of Lemma \ref{algebra nonzero} which asserts that the
$\pi$-torsion submodule of the dual Selmer group is trivial.

We should remark that one can prove the main results in this article
under the assumption that the dual Selmer group in question has the
property that it does not contain any nontrivial pseudo-null
submodule. For the case of dual Selmer groups of elliptic curves,
this property has been verified in many cases (see \cite{HO,
OcV02}). However, for dual Selmer groups of other Galois
representations, such pseudo-nullity results are less complete (but
see \cite{Oc, Su}). However, if one is willing to work with $p$-adic
extensions that are not totally real, we will see in this paper that
we do not require this property. It is also worthwhile mentioning
that we do not make use of any explicit formulas for the Akashi
series and Euler characteristics, and so our approaches differ from
those in \cite{Sh, SS14}.

We now give a brief description of the layout of the paper. In
Section \ref{algebra}, we recall certain algebraic notion which will
be used in the subsequent of the paper. In Section \ref{cyclotomic
section}, we review the main result of Matsuno \cite{Mat} on the
structure of the Selmer groups over the cyclotomic $\Zp$-extension.
In Section \ref{admissible section}, we review some of the results
in \cite{LimMHG} which we require. We will also show how Matsuno's
result can be applied to establish certain result on the structure
of the Selmer groups which is crucial in the proof of our main
results. In Section \ref{Akashi section}, we will prove our main
results on the Akashi series and Euler characteristics. In the final
section, we will prove the analogous result for the characteristic
elements of the Selmer groups.

\section{Algebraic Preliminaries} \label{algebra}

In this section, we recall some algebraic preliminaries that will be
required in the later part of the paper. Let $p$ be an odd prime.
Let $\Op$ be the ring of integers of a finite extension of $\Qp$.
Let $G$ be a compact $p$-adic Lie group without $p$-torsion. The
completed group algebra of $G$ over $\Op$ is given by
 \[ \Op\ps{G} = \plim_U \Op[G/U], \]
where $U$ runs over the open normal subgroups of $G$ and the inverse
limit is taken with respect to the canonical projection maps. It is
well known that $\Op\ps{G}$ is an Auslander regular ring (cf.
\cite[Theorem 3.26]{V02}). If $G$ is pro-$p$, then $\Op\ps{G}$ has
no zero divisors (cf.\ \cite{Neu}). Therefore, $\Op\ps{G}$ admits a
skew field $K(G)$ which is flat over $\Op\ps{G}$ (see \cite[Chapters
6 and 10]{GW} or \cite[Chapter 4, \S 9 and \S 10]{Lam}). If $M$ is a
finitely generated $\Op\ps{G}$-module, we define the
$\Op\ps{G}$-rank of $M$ to be
$$ \rank_{\Op\ps{G}}(M)  = \dim_{K(G)} \big(K(G)\ot_{\Op\ps{G}}M\big). $$
 We say that the module $M$ is a
\textit{torsion} $\Op\ps{G}$-module if $\rank_{\Op\ps{G}} M = 0$.

We continue to assume that $G$ is pro-$p$. Fix a local parameter
$\pi$ for $\Op$ and denote the residue field of $\Op$ by $k$. The
completed group algebra of $G$ over $k$ is defined similarly as
above. For a $k\ps{G}$-module $N$, we then define its $k\ps{G}$-rank
by
 $$ \rank_{k\ps{G}}(N)=
 \displaystyle\frac{\rank_{k\ps{G_0}}(N)}{|G:G_0|}, $$
where $G_0$ is an open normal uniform pro-$p$ subgroup of $G$. This
is integral and independent of the choice of $G_0$ (see
\cite[Proposition 1.6]{Ho}). Similarly, we will say that that $N$ is
a \textit{torsion} $\mathbb{F}_p\ps{G}$-module if
$\rank_{\mathbb{F}_p\ps{G}}N = 0$.

For a given finitely generated $\Op\ps{G}$-module $M$, we denote
$M(\pi)$ to be the $\Op\ps{G}$-submodule of $M$ consisting of
elements of $M$ which are annihilated by some power of $\pi$. Since
the ring $\Op\ps{G}$ is Noetherian, the module $M(\pi)$ is finitely
generated over $\Op\ps{G}$. Therefore, one can find an integer
$r\geq 0$ such that $\pi^r$ annihilates $M(\pi)$. Following
\cite[Formula (33)]{Ho}, we define
  \[\mu_{\Op\ps{G}}(M) = \sum_{i\geq 0}\rank_{k\ps{G}}\big(\pi^i
   M(\pi)/\pi^{i+1}\big). \]
(For another alternative, but equivalent, definition, see
\cite[Definition 3.32]{V02}.) By the above discussion, the sum on
the right is a finite one.

We introduce another algebraic invariant which was first defined in
\cite{CSS}. As before, $G$ is a compact (not necessarily pro-$p$)
$p$-adic Lie group without $p$-torsion. Let $H$ be a closed normal
subgroup of $G$ with $\Ga:= G/H\cong \Zp$. We say that the
\textit{Akashi series} of $M$ exists if $H_i(H,M)$ is
$\Op\ps{\Ga}$-torsion for every $i$. In the case of this event, we
denote $Ak_H(M)$ to be the \emph{Akashi series} of $M$ which is
defined to be \[ \displaystyle \prod_{i\geq 0} g_i^{(-1)^i},\]
 where $g_i$ is the characteristic polynomial of $H_i(H,M)$.
 Of course, the Akashi series is only well-defined up to a
unit in $\Op\ps{\Ga}$. Also, note that since $G$ (and hence $H$) has
no $p$-torsion, $H$ has finite $p$-cohomological dimension, and
therefore, the alternating product is a finite product. For the
remainder of the section, identify $\Op\ps{\Ga}\cong \Op\ps{T}$
under a choice of a generator of $\Ga$. A polynomial $T^n +
c_{n-1}T^{n-1} + \cdots +c_0$ in $\Op[T]$ is said to be a
\textit{Weierstrass polynomial} if $\pi$ divides $c_i$ for every
$0\leq i \leq n-1$. We now record the following proposition which
gives a relation between the Akashi series and various Iwasawa
invariants when $G$ is pro-$p$.

\bl \label{akashi mu} Let $G$ be a compact pro-$p$ $p$-adic group
without $p$-torsion, and let $H$ be a closed normal subgroup of $G$
with $\Ga:=G/H\cong \Zp$. Let $M$ be a finitely generated torsion
$\Op\ps{G}$-module which satisfies the following properties:
\begin{enumerate}
\item[$(i)$] $H_0(H,M)$ is a finitely generated torsion
$\Op\ps{G}$-module.
\item[$(ii)$] $H_i(H,M)$ is finitely generated over
$\Op$ for every $i\geq 1$.
\end{enumerate}
Then we have
 \[ Ak_H(M) = \pi^{\mu_{\Op\ps{\Ga}}(H_0(H,M)))}\frac{f}{g}, \]
 where $f$ and $g$ are Weierstrass polynomials. \el

\bpf
 This is an easy exercise which we leave to the reader.
\epf

We will also record the following proposition which is a special
case of \cite[Proposition 5.4]{LimMHG}.

\bp \label{algebra prop} Let $G$ be a compact pro-$p$ $p$-adic group
without $p$-torsion, and let $H$ be a closed normal subgroup of $G$
with $\Ga:= G/H\cong \Zp$. Let $M$ be a finitely generated
$\Op\ps{G}$-module which satisfies all of the following properties.
\begin{enumerate}
\item[$(i)$] $H_0(H,M)$ is a finitely generated torsion
$\Op\ps{G}$-module.
\item[$(ii)$] $H_i(H,M)$ is finitely generated over
$\Op$ for every $i\geq 1$.
\item[$(iii)$] $Ak_H(M)$ lies in $\Op\ps{\Ga}^{\times}$.
\end{enumerate}
Then $M$ is a finitely generated torsion $\Op\ps{H}$-module.   \ep

We say that the $G$-\textit{Euler characteristics} of an
$\Op\ps{G}$-module $M$ exists if $H_i(G,M)$ is finite for each
$i\geq 0$. In the event of such, the $G$-\textit{Euler
characteristics} is given by
\[ \chi(G,M) = \prod_{i\geq 0}|H_i(G,M)|^{(-1)^i}. \]
Again, since $G$ has no $p$-torsion, the above product is a finite
one. The connection between the Akashi series and the Euler
characteristics is given as follow.

\bp \label{Akashi Euler}
 Let $G$ be a compact $p$-adic group
without $p$-torsion, and let $H$ be a closed normal subgroup of $G$
with $G/H\cong \Zp$. Let $M$ be a finitely generated
$\Op\ps{G}$-module whose $G$-Euler characteristics is well defined.
Then the Akashi series of $M$ exists and we have
 \[  \chi(G,M) = |\varphi(Ak_H(M))|_{\pi}^{-1}, \]
 where $|~|_{\pi}$ is the $\pi$-adic norm with $|\pi|_{\pi} = q^{-1}$ $($here $q$ is
 the order of $k)$,
 and $\varphi$ is the augmentation map from $\Op\ps{\Ga}$ to $\Op$.
\ep

\bpf
 See \cite[Theorem 3.6]{CFKSV} or \cite[Lemma 4.2]{CSS}.
\epf

\section{Selmer groups over cyclotomic $\Zp$-extension}
\label{cyclotomic section}

As before, $p$ will denote an odd prime and $\Op$ will denote the
ring of integers of a fixed finite extension $K$ of $\Qp$. Suppose
that we are given the following datum $\big(A, \{A_v\}_{v|p})$
defined over a number field $F$:

\begin{enumerate}
 \item[(C1)] $A$ is a
cofinitely generated cofree $\Op$-module of $\Op$-corank $d$ with a
continuous, $\Op$-linear $\Gal(\bar{F}/F)$-action which is
unramified outside a finite set of primes of $F$.

 \item[(C2)] For each prime $v$ of $F$ above $p$, $A_v$ is a
$\Gal(\bar{F}_v/F_v)$-submodule of $A$ which is cofree of
$\Op$-corank $d_v$.

\item[(C3)] For each real prime $v$ of $F$, we write $A_v^+=
A^{\Gal(\bar{F}_v/F_v)}$  which is assumed to be cofree of
$\Op$-corank $d^+_v$.

\item[(C4)] The following equality
  \begin{equation} \label{data equality}
  \sum_{v|p} (d-d_v)[F_v:\Qp] = dr_2(F) +
 \sum_{v~\mathrm{real}}(d-d^+_v)
  \end{equation}
holds. Here $r_2(F)$ denotes the number of complex primes of $F$.
\end{enumerate}

We now consider the base change property of our datum. Let $L$ be a
finite extension of $F$. We can then obtain another set of datum
$\big(A, \{A_w\}_{w|p}\big)$ over $L$ as follows: we consider $A$ as
a $\Gal(\bar{F}/L)$-module, and for each prime $w$ of $L$ above $p$,
we set $A_w =A_v$, where $v$ is a prime of $F$ below $w$, and view
it as a $\Gal(\bar{F}_v/L_w)$-module. Then $d_w = d_v$. For each
real prime $w$ of $L$, one sets $A^{\Gal(\bar{L}_w/L_w)}=
A^{\Gal(\bar{F}_v/F_v)}$ and writes $d^+_w = d^+_v$, where $v$ is a
real prime of $F$ below $w$. In general, the $d_w$'s and $d_w^+$'s
need not satisfy equality (C4). We now record the following lemma
which gives some sufficient conditions for the equality in (C4) to
hold for the datum $\big(A, \{A_w\}_{w|p}\big)$ over $L$ (see
\cite[Lemma 4.1]{LimMu} for the proofs).

\bl \label{data base change} Suppose that $\big(A,
\{A_v\}_{v|p}\big)$ is a datum defined over $F$ which satisfies
$(C1)-(C4)$.
 Suppose further that at least one of the following statements holds.
 \begin{enumerate}
\item[$(i)$] All the archimedean primes of $F$ are unramified in
$L$.

\item[$(ii)$] $[L:F]$ is odd

\item[$(iii)$] $F$ is totally imaginary.

\item[$(iv)$] $F$ is totally real, $L$ is totally imaginary and
\[ \sum_{v~\mathrm{real}} d^+_v = d[F:Q]/2.\]
 \end{enumerate}
Then we have the equality
 \[ \sum_{w|p} (d-d_w)[L_w:\Qp] = dr_2(L) +
 \sum_{w~\mathrm{real}}(d-d^+_w).\]
 \el

We now introduce the Selmer groups. Let $S$ be a finite set of
primes of $F$ which contains all the primes above $p$, the ramified
primes of $A$ and all infinite primes. Denote $F_S$ to be the
maximal algebraic extension of $F$ unramified outside $S$ and write
$G_S(\mathcal{L}) = \Gal(F_S/\mathcal{L})$ for every algebraic
extension $\mathcal{L}$ of $F$ which is contained in $F_S$. Let $L$
be a finite extension of $F$ contained in $F_S$ such that the datum
$\big(A, \{A_w\}_{w|p}\big)$ satisfies (C1)$-$(C4). For a prime $w$
of $L$ lying over $S$, set
\[ H^1_{str}(L_w, A)=
\begin{cases} \ker\big(H^1(L_w, A)\lra H^1(L_w, A/A_w)\big) & \text{\mbox{if} $w$
 divides $p$},\\
 \ker\big(H^1(L_w, A)\lra H^1(L^{ur}_w, A)\big) & \text{\mbox{if} $w$ does not divide $p$,}
\end{cases} \]
 where $L_w^{ur}$ is the maximal unramified extension of $L_w$.
The (strict) Selmer group attached to the datum is then defined by
\[ \Sel^{str}(A/L) := \ker\Bigg( H^1(G_S(L),A)\lra
\bigoplus_{w \in S_L}\frac{H^1(L_w, A)}{H^1_{str}(L_w,
A)_{div}}\Bigg).\] Here $N_{div}$ denotes the maximal
$\Op$-divisible submodule of $N$ and $S_L$ denotes the set of primes
of $L$ above $S$. We now denote $F^{\cyc}$ to be the cyclotomic
$\Zp$-extension of $F$, and $F_n$ its intermediate subextension with
$[F_n:F]=p^n$. Write $\Ga = \Gal(F_{\infty}/F)$ and $\Ga_n
=\Gal(F_{\infty}/F_n)$. We define the Selmer group over $F^{\cyc}$
by
\[\Sel^{str}(A/F^{\cyc})=
\ilim_n \Sel^{str}(A/F_n),\] where the limit runs over the
intermediate extensions $F_n$ of $F$ contained in $F^{\cyc}$. We
will write $X(A/F^{\cyc})$ for the Pontryagin dual of
$\Sel^{str}(A/F^{\cyc})$.

We now introduce two additional conditions on our datum.

(Fin$_p$) For every prime $w$ of $F^{\cyc}$ above $p$,
$H^0(F_w^{\cyc}, A_v)$ is finite. Here $v$ is a prime of $F$ under
$w$.

(Fin) $H^0(F^{\cyc}, A^*)$ is finite. Here $A^* =
\Hom_{\cts}(T_{\pi}(A),\mu_{p^{\infty}})$, where $T_{\pi}(A) =
\plim_i A[\pi^i]$.

The following lemma gives an alternative description of the Selmer
group over $F^{\cyc}$.

\bl \label{Selmer limit cyc}
 Suppose that $(Fin_p)$ holds. Then we have
 \[ \Sel^{str}(A/F^{\cyc}) = \ker\Big(H^1(G_S(F^{\cyc}), A)
 \lra \bigoplus_{v\in S} J_v(A/F^{\cyc})
 \Big), \]
 where
 \[ J_v(A/F^{\cyc}) = \begin{cases}
 \bigoplus_{w|v}H^1(F^{\cyc}_w, A/A_v),& \mbox{if } v\mbox{ divides }p \\
   \bigoplus_{w|v}H^1(F^{\cyc}_w, A), & \mbox{if } v\mbox{ does not divides }p \end{cases}
 \] \el

\bpf
 Since (Fin$_p$) holds, it follows from the discussion in
 \cite[P 111]{G89} that $H^1(F^{\cyc}_w, A_v)$ is divisible for
 $w|p$. One can then proceed using a similar argument to
 that in \cite[Lemma 3.2]{Mat}.
\epf

For the remainder of the section, we will review the result of
Matsuno on the structure of the dual Selmer group which will play an
important role in proving the main results of the paper. Now
following \cite{Mat}, we define the following finite
$\Op\ps{\Ga}$-module.

\bd For each prime $w$ of $F^{\cyc}$, we define a subgroup of
$H^0(F^{\cyc}_w, A^*)$ by
\[
 B_w= \begin{cases}
 H^0(F^{\cyc}_w, A^*_v),& \mbox{if } w \mbox{ lies above } v|p, \\
   H^0(F^{\cyc}_w, A^*_v), & \mbox{if } w \mbox{ is archimedean, } \\
 H^0(F^{\cyc}_w, A^*_v)_{div}, & \mbox{otherwise. }  \end{cases}
 \]
Then we define a subgroup $B(F^{\cyc}, A^*)$ of $H^0(F^{\cyc}, A^*)$
by
\[ B(F^{\cyc}, A^*)= \bigcap_w i_w^{-1}(B_w),\]
 where $w$ runs over all primes of $F^{\cyc}$ and $i_w$ is the
 canonical injection
 \[ H^0(F^{\cyc}, A^*)\hookrightarrow H^0(F^{\cyc}_w, A^*). \]
\ed

For a finitely generated torsion $\Op\ps{\Ga}$-module $M$, we denote
$Fin_{\Op\ps{\Ga}}(M)$ to be the maximal finite submodule of $M$.
(This is well-defined; for instance, see \cite[Lemma 3.2]{Ch}.) The
following lemma gives a relationship between $B(F^{\cyc}, A^*)$ and
the maximal finite submodule of $X(A/F^{\cyc})$.

\bl
 Suppose that $(Fin_p)$ and $(Fin)$ hold. Assume that $X(A/F^{\cyc})$ is a
 torsion $\Op\ps{G}$-module. Then there is
 an $\Op\ps{\Ga}$ injection
 \[ Fin_{\Op\ps{\Ga}} X(A/F^{\cyc}) \hookrightarrow B(F^{\cyc}, A^*).\] \el

\bpf
 See \cite[Proposition 3.9]{Ch} or \cite[Corollary 4.3]{Mat}.
\epf

\bp
 Let $p$ be an odd prime. Assume that $(Fin_p)$ and $(Fin)$ hold.
 Suppose that $F$ is not totally real and that $X(A/F^{\cyc})$
 has a nontrivial finite $\Op\ps{\Ga}$-module. Then $X(A/F^{\cyc})$
 is not finitely generated over $\Op$.
\ep

\bpf
 This is proven by a similar argument to that in
 \cite[Theorem 3.11]{Ch} or  \cite[Proposition
 7.5]{Mat} by appealing to the preceding lemma
 and the fact that the Galois group of the maximal abelian pro-$p$
 extension of $F^{\cyc}$ unramified outside $p$ has positive $\Zp\ps{\Ga}$-rank
 when $F$ is not totally real. \epf

We record a corollary which will play a crucial role in the
subsequent of the paper.

\bc \label{finite coro}
 Let $p$ be an odd prime. Assume that $(Fin_p)$ and $(Fin)$ hold.
 Suppose that $F$ is not totally real and that
 $X(A/F^{\cyc})$ is finitely generated over $\Op$. Then $X(A/F^{\cyc})$
 is $\Op$-torsionfree. \ec

We end the section mentioning two basic examples of our datum.

(i) $A = \mathcal{A}[p^{\infty}]$, where $\mathcal{A}$ is an abelian
variety defined over an arbitrary finite extension $F$ of $\Q$ with
good ordinary reduction at all places $v$ of $F$ dividing $p$. For
each $v|p$, it follows from \cite[P. 150-151]{CG} that we have a
$\Gal(\bar{F}_v/F_v)$-submodule $A_v$ which can be characterized by
the property that $A/A_v$ is the maximal
$\Gal(\bar{F}_v/F_v)$-quotient of $\mathcal{A}[p^{\infty}]$ on which
some subgroup of finite index in the inertia group $I_v$ acts
trivially. It is not difficult to verify that (C1)-(C4) are
satisfied. The condition (Fin$_p$) is also known to be satisfied
(see \cite[Section 5]{Mat}) and the condition (Fin) is a well-known
consequence of a theorem of Imai \cite{Im}.

(ii) Let $V$ be the Galois representation attached to a primitive
Hecke eigenform $f$ for $GL_2 /\Q$, which is ordinary at $p$,
relative to some fixed embedding of the algebraic closure of $\Q$
into $\Qp$. By the work of Mazur-Wiles \cite{MW}, $V$ contains a
one-dimensional $\Qp$-subspace $V_v$ invariant under
$\Gal(\bar{\Q}_p/\Qp)$ with the property that the inertial subgroup
$I_p$ acts via a power of the cyclotomic character on $V_v$ and
trivially on $V/V_v$. By compactness, $V$ will always contain a free
$\Op$-submodule $T$, which is stable under the action of
$\Gal(\bar{F} /F)$. For such an $\Op$-lattice $T$, we write $A =
V/T$ and $A_v = V_v/ (T\cap V_v)$. The condition (Fin$_p$) can be
easily seen to follow from the property that the inertial subgroup
$I_p$ acts via a power of the cyclotomic character on $V_v$ (or see
the proof of \cite[Lemma 3.6]{Sh}). The condition (Fin) is shown in
the proof of \cite[Lemma 2.2]{Su}.

\section{Selmer groups over strongly admissible $p$-adic Lie extensions}
\label{admissible section}

We say that $F_{\infty}$ is a \textit{strongly admissible $p$-adic
Lie extension} of $F$ if (i) $\Gal(F_{\infty}/F)$ is a compact
$p$-adic Lie group without $p$-torsion, (ii) $F_{\infty}$ contains
the cyclotomic $\Zp$-extension $F^{\cyc}$ of $F$ and (iii)
$F_{\infty}$ is unramified outside a finite set of primes of $F$. We
will write $G = \Gal(F_{\infty}/F)$, $H = \Gal(F_{\infty}/F^{\cyc})$
and $\Ga =\Gal(F^{\cyc}/F)$.

Let $\big(A, \{A_v\}_{v|p}\big)$ be a datum defined as in the
previous section which satisfies (C1)$-$(C4). \textit{For the
remainder of the paper, we will always assume that for every finite
extension $L$ of $F$ in $F_{\infty}$, the datum $\big(A,
\{A_w\}_{w|p}\big)$ over $L$ obtained by base change satisfies}
(C1)$-$(C4). Let $S$ be a finite set of primes which contains all
the primes above $p$, the ramified primes of $A$, the archimedean
primes and the primes that are ramified in $F_{\infty}/F$. The
Selmer group of the data over $F_{\infty}$ is defined by
\[ \Sel^{str}(A/F_{\infty}) = \ilim_L \Sel^{str}(A/L), \]
where the limit runs over all the finite extensions $L$ of $F$
contained in $F_{\infty}$. By a standard argument as in
\cite[Section 2]{CS12}, one can show that the Selmer group is
independent of the choice of $S$ as long as it contains all the
primes above $p$, the ramified primes of $A$, the archimedean primes
and the primes that are ramified in $F_{\infty}/F$. We write
$X(A/F_{\infty})$ for the Pontryagin dual of
$\Sel^{str}(A/F_{\infty})$. We record the following proposition.

\bp \label{Selmer limit} Suppose that $(Fin_p)$ holds for every
finite extension $L$ of $F$ contained in $F_{\infty}$. Then we have
the following description of the strict Selmer groups.
\[ \Sel^{str}(A/F_{\infty}) = \ker\Big(H^1(G_S(F_{\infty}), A)
 \stackrel{\la_{A/F_{\infty}}}{\lra} \bigoplus_{v\in S} J_v(A/F_{\infty})
 \Big), \]
 where
 \[ J_v(A/F_{\infty}) = \ilim_L J_v(A/L^{\cyc}). \]
Here the limit is taken over all the finite extensions $L$ of $F$
contained in $F_{\infty}$.  \ep

\bpf
 For each $L$, we denote $L_n$ to be the intermediate extension of
 $L$ contained in $L^{\cyc}$ with $[L:L_n]=p^n$. Then we have
 \[ \Sel^{str}(A/F_{\infty}) = \ilim_L\ilim_n \Sel^{str}(A/L_n) =
 \ilim_L \Sel^{str}(A/L^{\cyc}).  \]
 The conclusion of the proposition is now immediate from
Lemma \ref{Selmer limit cyc}. \epf

For the remainder of the paper, we will also impose the following
condition.

\medskip
\noindent \textbf{(Fin$_{\infty}$)}: The data $\big(A,
\{A_v\}_{v|p}\big)$ satisfies (Fin$_p$) and (Fin) for every finite
extension $L$ of $F$ contained in $F_{\infty}$.

\bl \label{algebra nonzero}  Assume that $\mathbf{(Fin_{\infty})}$
is valid, and assume that $F_{\infty}$ is not totally real. Then the
following statements hold.
\begin{enumerate}
\item[$(a)$] If $X(A/F_{\infty})\neq 0$, then $X(A/F_{\infty})$
is not a finitely generated torsion $\Op\ps{H}$-module.
\item[$(b)$] If $X(A/F_{\infty})$ is finitely generated over
$\Op\ps{H}$, then $X(A/F_{\infty})(\pi) =0$.
 \end{enumerate}
  \el
\bpf (a) This is proven by a similar argument to that in \cite[Lemma
5.8]{LimMHG} by appealing to Corollary \ref{finite coro}.

(b) Suppose that $X(A/F_{\infty})$ is finitely generated over
$\Op\ps{H}$. In particular, this implies that $X(A/L^{\cyc})$ is
finitely generated over $\Op$ for every finite extension $L$ of $F$
contained in $F_{\infty}$ such that $\Gal(F_{\infty}/L)$ is pro-$p$.
Since $F_{\infty}$ is not totally real, we can rewrite
$X(A/F_{\infty}) = \plim_LX(A/L^{\cyc})$, where $L$ runs through all
finite extensions of $F$ contained in $F_{\infty}$ such that
$\Gal(F_{\infty}/L)$ is pro-$p$ and $L$ is not totally real. By
Corollary \ref{finite coro}, the multiplication by $\pi$-map on
$X(A/L^{\cyc})$ is an injection. Since inverse limit is left exact,
it follows that the multiplication by $\pi$-map on $X(A/F_{\infty})$
is also an injection.
 \epf

 We will write
$\la_{A/F_{\infty}}$ for the localization map
\[ H^1(G_S(F_{\infty}), A)
 \lra \bigoplus_{v\in S} J_v(A/F_{\infty}) \] as given in Proposition
 \ref{Selmer limit}. We can now record the following result.

\bp \label{Selmer fg} Assume that $(i)$ $\mathbf{(Fin_{\infty})}$
holds, $(ii)$ $G=\Gal(F_{\infty}/F)$ is pro-$p$ and has no
$p$-torsion, $(iii)$ $X(A/F^{\cyc})$ is $\Op\ps{\Ga}$-torsion,
$(iv)$ $H^2(G_S(F_{\infty}),A)=0$ and $(v)$ $\la_{A/F_{\infty}}$ is
surjective. Then the following statements hold.

\begin{enumerate}
\item[$(a)$] $H_0(H, X(A/F_{\infty}))$ is $\Op\ps{\Ga}$-torsion and
its $\mu_{\Op\ps{\Ga}}$-invariant is precisely
$\mu_{\Op\ps{\Ga}}\big(X(A/F^{\cyc})\big)$.

\item[$(b)$] $H_i(H, X(A/F_{\infty}))$ is finitely generated over
$\Op$ for every $i\geq 1$.
\item[$(c)$] $H_i(H, X(A/F_{\infty}))=0$  for $i\geq \dim H$.
\end{enumerate}

In particular, we have
\[ Ak_H\big(X(A/F_{\infty})\big) =
\pi^{\mu_{\Op\ps{\Ga}}\big(X(A/F^{\cyc})\big)}\frac{f}{g}\]
 for some Weierstrass polynomials $f$ and $g$. \ep

\bpf
 The assertions on the $H$-homology of $X(A/F_{\infty})$ can
be proven by a similar argument to that in \cite[Proposition
5.1]{LimMHG}. The assertion on the Akashi series is then an
immediate consequence of this and Lemma \ref{akashi mu}. \epf

\br
 It follows from a standard argument that
conditions (iv) and (v) will imply that $X(A/F_{\infty})$ is a
torsion $\Op\ps{G}$-module. Also, there are many cases of conditions
(iv) and (v) following from (iii) and we mention a few of them here.

(a) Suppose $G$ has dimension 2. By our (Fin) assumption, we may
apply a standard argument (for instance, see \cite[Proposition
3.3(a)]{LimMHG}) to show that $H^2(G_S(F^{\cyc}),A)=0$ and
$\la_{A/F^{\cyc}}$ is surjective. One can now apply an argument
similar to \cite[Theorem 7.4]{SS} to obtain the conclusion that
$H^2(G_S(F_{\infty}),A)=0$ and $\la_{A/F_{\infty}}$ is surjective.

(b) $G$ is a solvable uniform pro-$p$ group and $A(F_{\infty})$ is
finite. Then by a similar argument to that in \cite[Theorem
2.3]{HO}, we have that $X(A/F_{\infty})$ is a torsion
$\Op\ps{G}$-module. By a standard argument (for instances, see
\cite[Proposition 3.3(a)]{LimMHG}), we obtain the conclusion that
$H^2(G_S(F_{\infty}),A)=0$ and $\la_{A/F_{\infty}}$ is surjective.

(c) $G$ is a solvable uniform pro-$p$ group, and for each $v\in S$,
the decomposition group of $G = \Gal(F_{\infty}/F)$ at $v$ has
dimension $\geq 2$. Again by the argument in \cite[Theorem 2.3]{HO},
we have that $X(A/F_{\infty})$ is a torsion $\Op\ps{G}$-module. The
required conclusion will now follow from the argument in
\cite[Proposition 3.3(b)]{LimMHG}.

(d) $X(A/F^{\cyc})$ is finitely generated over $\Op$. Then for every
finite extension $L$ of $F$ in $F_{\infty}$, it follows from the
argument in \cite{HM} that $X(A/L^{\cyc})$ is finitely generated
over $\Op$, and in particularly, torsion over $\Op\ps{\Ga_L}$, where
$\Ga_L = \Gal(L^{\cyc}/L)$. Combining this with our (Fin)
assumption, we can apply a similar argument to that in
\cite[Proposition 3.3(i), Corollary 3.4]{LimMHG} to conclude that
$H^2(G_S(F_{\infty}),A)=0$ and $\la_{A/F_{\infty}}$ is surjective.
\er

\section{Comparing Akashi series of Selmer groups} \label{Akashi
section}

Retain the notation from the previous section. To state the main
result of this section, we need to introduce another datum $\big(B,
\{B_v\}_{v|p}\big)$, where we assume to satisfy conditions
(C1)-(C4). We will also assume that the datum obtained by base
change over every finite extension of $L$ of $F$ in $F_{\infty}$
also satisfy conditions (C1)-(C4). From now on, $S$ will always
denote a finite set of primes which contains all the primes above
$p$, the ramified primes of $A$ and $B$, the ramified primes of
$F_{\infty}/F$ and the archimedean primes. We introduce the
following important congruence condition on $A$ and $B$ which allows
us to be able to compare the Selmer groups of $A$ and $B$.

\smallskip \noindent $\mathbf{(Cong)}$ : There is an isomorphism
$A[\pi]\cong B[\pi]$ of $G_S(F)$-modules which induces a
$\Gal(\bar{F}_v/F_v)$-isomorphism $A_v[\pi]\cong B_v[\pi]$ for every
$v|p$.

\medskip
We introduce some notation which we need in our discussion.
 We write $C_v =A$ for $v\nmid p$ and $C_v =
A/A_v$ for $v|p$. Similarly, we write $D_v =B$ for $v\nmid p$ and
$D_v = B/B_v$ for $v|p$. Let $S_2$ denote the set of primes in $S$
such that the decomposition group of $G= \Gal(F_{\infty}/F)$ at $v$
has dimension $\geq 2$. For every prime $w$ of $F^{\cyc}$ above $v$,
we write $C_v(F^{\cyc}_w)= (C_v)^{\Gal(\bar{F}_v/F^{\cyc}_w)}$ and
$D_v(F^{\cyc}_w)= (D_v)^{\Gal(\bar{F}_v/F^{\cyc}_w)}$.

\bt \label{Akashi congruent}
 Let $F_{\infty}$ be a strongly admissible pro-$p$ $p$-adic Lie extension of $F$
 of dimension at least 2.
 Assume that $F$ is not totally real.
 Suppose that all the hypotheses in Proposition \ref{Selmer fg} are valid
 for $A$ and $B$, and that $\mathbf{(Cong)}$ is satisfied. Furthermore,
 suppose that for every $v\notin S_2$, we have
 \[\dim_k \big(C_v(F^{\cyc}_w)/\pi\big) = \dim_k\big(D_v(F^{\cyc}_w)/\pi\big)\]
  for every prime $w$ of $F^{\cyc}$ above $v$.

Then $Ak_H(X(A/F_{\infty}))$ is a unit in $\Op\ps{\Ga}$ if and only
if $Ak_H(X(B/F_{\infty}))$ is a unit in $\Op\ps{\Ga}$. \et

\br
 We mention some cases when the final hypothesis of the theorem is
 satisfied.

 (a) $C_v(F^{\cyc}_w)$ and $D_v(F^{\cyc}_w)$ are
 $\Op$-divisible. This is immediate. We also note that
 this has been known to hold in certain cases
(see \cite[Lemma 4.1.3]{EPW}; also see the proofs of \cite[Theorem
4.7]{Sh} and \cite[Theorem 3.4]{SS14}).

 (b) $C_v(F^{\cyc}_w)$ and $D_v(F^{\cyc}_w)$ are finite. Then by Howson's formula
 \cite[Corollary 1.10]{Ho}, we have
 \[\dim_k \big(C_v(F^{\cyc}_w)/\pi\big) = \dim_k\big(C_v(F^{\cyc}_w)[\pi]\big).\]
Similarly, we have an analogue formula for $D$. By the
$\mathbf{(Cong)}$ assumption, we have
 \[\dim_k \big(C_v(F^{\cyc}_w)[\pi]\big) = \dim_k\big(D_v(F^{\cyc}_w)[\pi]\big).\]
 Combining all these equalities, we obtain the required conclusion.
 \er

In preparation of the proof of the theorem, we record a proposition.

\bp \label{Akashi congruent lemma}
 Retain the assumptions of Theorem \ref{Akashi congruent}, then
 $X(A/F_{\infty})$ is finitely generated over $\Op\ps{H}$ if and only
if $X(B/F_{\infty})$ is finitely generated over $\Op\ps{H}$.
Furthermore, we have
 \[ \rank_{\Op\ps{H}}\big(X(A/F_{\infty})\big) =
 \rank_{\Op\ps{H}}\big(X(B/F_{\infty})\big). \]\ep

We now prove Theorem \ref{Akashi congruent} assuming the validity of
Proposition \ref{Akashi congruent lemma}.

\bpf[Proof of Theorem \ref{Akashi congruent}]
 Suppose that $Ak_H(X(A/F_{\infty}))$ is a unit in $\Op\ps{\Ga}$.
It then follows from Propositions \ref{algebra prop} and \ref{Selmer
fg} that $X(A/F_{\infty})$ is a finitely generated torsion
$\Op\ps{H}$-module. By Proposition \ref{Akashi congruent lemma},
this in turn implies that $X(B/F_{\infty})$ is a finitely generated
torsion $\Op\ps{H}$-module. Since $F$ is not totally real, it
follows from Lemma \ref{algebra nonzero}(a) that $X(B/F_{\infty})
=0$. In particular, this implies that $Ak_H(X(B/F_{\infty}))$ is a
unit in $\Op\ps{\Ga}$. \epf

The remainder of the section will be devoted to the proof of
Proposition \ref{Akashi congruent lemma}. In preparation of the
proof, we first introduce the ``mod $\pi$" Selmer group.  For every
finite extension $\mathcal{F}$ of $F^{\cyc}$, we define
$J_v(A[\pi]/\mathcal{F})$ to be
 \[  \bigoplus_{w|v}H^1(\mathcal{F}_w, A[\pi])~\mbox{or}~
  \bigoplus_{w|v}H^1(\mathcal{F}_w, A/A_v[\pi])\]
according as $v$ does not or does divide $p$. We then define
\[ J_v(A[\pi]/F_{\infty}) = \ilim_{\mathcal{F}} J_v(A[\pi]/\mathcal{F}),\]
where the direct limit is taken over all finite extensions
$\mathcal{F}$ of $F^{\cyc}$ contained in $F_{\infty}$. The mod $\pi$
Selmer group is then defined by
 \[ \Sel^{str}(A[\pi]/F_{\infty}) = \ker\Big(H^1(G_S(F_{\infty}), A[\pi])
 \lra \bigoplus_{v\in S} J_v(A[\pi]/F_{\infty})
 \Big).\]
We write $X(A[\pi]/F_{\infty})$ for the Pontryagin dual of
$\Sel^{str}(A[\pi]/F_{\infty})$. We are now in the position to prove
Proposition \ref{Akashi congruent lemma}.

\bpf[Proof of Proposition \ref{Akashi congruent lemma}]
 Consider the
following diagram
\[  \entrymodifiers={!! <0pt, .8ex>+} \SelectTips{eu}{}\xymatrix{
    0 \ar[r]^{} & \Sel^{str}(A[\pi]/F_{\infty}) \ar[d]_{a} \ar[r] &
    H^1(G_S(F_{\infty}), A[\pi])
    \ar[d]_{b}
    \ar[r]^{\psi_{A}} & \bigoplus_{v\in S}J_v(A[\pi]/F_{\infty}) \ar[d]_{c} \\
    0 \ar[r]^{} & \Sel^{str}(A/F_{\infty})[\pi] \ar[r]^{}
    & H^1(G_S(F_{\infty}), A)[\pi] \ar[r] & \
    \bigoplus_{_{v\in S}}J_v(A/F_{\infty})[\pi]  } \]
with exact rows. It follows from a standard argument that $\ker a$
and $\coker a$ are cofinitely generated over $k\ps{H}$. Hence it
follows that $X(A/F_{\infty})$ is finitely generated over
$\Op\ps{H}$ if and only if $X(A[\pi]/F_{\infty})$ is finitely
generated over $k\ps{H}$. One also has a similar assertion for $B$.
Now by the $\mathbf{(Cong)}$ assumption, we have
$X(A[\pi]/F_{\infty})\cong X(B[\pi]/F_{\infty})$. Therefore, the
first assertion of the proposition is established.

We proceed to show the second assertion. By Lemma \ref{algebra
nonzero}(b), we have $X(A/F_{\infty})[\pi]=0$, or equivalently,
$\Sel^{str}(A/F_{\infty})/\pi=0$. Therefore, the above diagram can
be completed to the following.

\[  \entrymodifiers={!! <0pt, .8ex>+} \SelectTips{eu}{}\xymatrix{
    0 \ar[r]^{} & \Sel^{str}(A[\pi]/F_{\infty}) \ar[d]_{a} \ar[r] &
    H^1(G_S(F_{\infty}), A[\pi])
    \ar[d]_{b}
    \ar[r]^{\psi_{A}} & \bigoplus_{v\in S}J_v(A[\pi]/F_{\infty}) \ar[d]_{c} \ar[r] & \coker\psi_{A} \ar[r] \ar[d] &0\\
    0 \ar[r]^{} & \Sel^{str}(A/F_{\infty})[\pi] \ar[r]^{}
    & H^1(G_S(F_{\infty}), A)[\pi] \ar[r] & \
    \bigoplus_{_{v\in S}}J_v(A/F_{\infty})[\pi] \ar[r] & 0 & } \]
 Since $\ker b$ and $\ker c$ are cofinitely generated over $k\ps{H}$,
and $b$ and $c$ are surjective, it follows from a simple diagram
chasing argument that $\coker\psi_A$ is cofinitely generated over
$k\ps{H}$. Furthermore, $\ker b \cong A(F_{\infty})/\pi$ is
cofinitely generated over $k$, and hence is a cotorsion
$k\ps{H}$-module. Therefore, by another diagram chasing argument, we
have
\[ \mathrm{corank}_{k\ps{H}}\big(\Sel^{str}(A/F_{\infty})[\pi]\big)
 =\mathrm{corank}_{k\ps{H}}(\ker c)
 + \mathrm{corank}_{k\ps{H}}\big(\Sel^{str}(A[\pi]/F_{\infty})\big) -
 \mathrm{corank}_{k\ps{H}}\big(\coker\psi_A\big).
\]
 Again, since $\Sel^{str}(A/F_{\infty})/\pi=0$, it follows from
\cite[Corollary 1.10]{Ho} that
\[ \mathrm{corank}_{\Op\ps{H}}\big(\Sel^{str}(A/F_{\infty})\big) =
\mathrm{corank}_{k\ps{H}}\big(\Sel^{str}(A/F_{\infty})[\pi]\big). \]

We now compute $\mathrm{corank}_{k\ps{H}}(\ker c)$. We write $c
=\oplus_w c_w$, where $w$ runs over the set of primes of $F^{\cyc}$
above $S$. Denote $H_w$ to be the decomposition group of
$F_{\infty}/F^{\cyc}$ corresponding to a fixed prime of $F_{\infty}$
above $w$. Then we have $\coker c_w =
\mathrm{Coind}^{H_w}_H\big(C_v(F_{\infty,w})/\pi\big)$, where $v$ is
the prime of $F$ below $w$.  Now it is easy to verify that
\[ \mathrm{corank}_{k\ps{H}}(\ker c_w)=
\begin{cases} \mathrm{rank}_{k}\big(C_v(F_{\infty,w})/\pi\big) &
 \text{if $v\notin S_2$},\\
0 & \text{if $v\in S_2$.}
\end{cases} \]
 In fact, since $F_{\infty}/F$ is a strongly admissible pro-$p$
 $p$-adic Lie extension, we have $F_{\infty,w} = F^{\cyc}_w$ when $v\notin S_2$. Therefore, in
 conclusion, we have
\[ \rank_{\Op\ps{H}}\big(X(A/F_{\infty})\big)
 = \sum_{v\notin S_2}\sum_{w|v}\dim_k \big(C_v(F^{\cyc}_w)/\pi\big)
 + \mathrm{corank}_{k\ps{H}}\big(\Sel^{str}(A[\pi]/F_{\infty})\big) -
 \mathrm{corank}_{k\ps{H}}\big(\coker\psi_A\big).
\]
By a similar argument, we have the same equality for $B$. By the
$\mathbf{(Cong)}$ assumption, we have
\[ \mathrm{corank}_{k\ps{H}}\big(\Sel^{str}(A[\pi]/F_{\infty})\big) -
 \mathrm{corank}_{k\ps{H}}\big(\coker\psi_A\big)
 = \mathrm{corank}_{k\ps{H}}\big(\Sel^{str}(B[\pi]/F_{\infty})\big) -
 \mathrm{corank}_{k\ps{H}}\big(\coker\psi_B\big).\] Combining this with the
 hypothesis of the proposition, we obtain
 \[ \rank_{\Op\ps{H}}\big(X(A/F_{\infty})\big) =
 \rank_{\Op\ps{H}}\big(X(B/F_{\infty})\big)\] as required. \epf

As a corollary, we have the analog result for the $G$-Euler
characteristics of the dual Selmer groups.

 \bc \label{Euler congruent}
 Retain the assumptions of Theorem \ref{Akashi congruent}. Assume further that
$H_i(H,X(A/F_{\infty}))$ and $H_i(H,X(B/F_{\infty}))$ are finite for
$i\geq 1$. Then we have $\chi(G, X(A/F_{\infty}))=1$ if and only if
$\chi(G, X(A/F_{\infty}))=1$. \ec

\bpf
 By the hypothesis that $H_i(H,X(A/F_{\infty}))$ are finite for
$i\geq 1$, we have that $Ak_H \big(X(A/F_{\infty})\big)$ is
precisely the characteristic polynomial of $H_0(H,
X(A/F_{\infty}))$. In particular, we have
$Ak_H\big(X(A/F_{\infty})\big) \in\Op\ps{\Ga}$. Combining this
observation with Proposition \ref{Akashi Euler}, we have $\chi(G,
X(A/F_{\infty}))=1$ if and only $Ak_H \big(X(A/F_{\infty})\big)$ is
a unit in $\Op\ps{\Ga}$. Similarly, we have $\chi(G,
X(B/F_{\infty}))=1$ if and only $Ak_H\big( X(B/F_{\infty})\big)$ is
a unit in $\Op\ps{\Ga}$. The conclusion of the theorem will now
follow from Theorem \ref{Akashi congruent}. \epf

\br
 We note that the finiteness condition on the $H$-homology in
 Corollary \ref{Euler congruent} is automatically satisfied
 if $G$ has dimension 2 (cf. Proposition \ref{Selmer
 fg}(c)). \er

We end the section proving a version of Theorem \ref{Akashi
congruent} for Selmer groups twisted by certain Artin
representations over a strongly admissible $p$-adic Lie extension
which is not pro-$p$. For the remainder of the section, we will
assume that our data $(A, \{A_v\}_{v|p})$ and $(B, \{B_v\}_{v|p})$
are defined over $\Q$. Let $F_{\infty}$ be a strongly admissible
pro-$p$ extension of $\Q(\mu_p)$ which contains
$\Q(\mu_{p^{\infty}}, m^{p^{-\infty}})$ for some $p$-power free
integer $m$. From now on, we write $G= \Gal(F_{\infty}/\Q)$,
$H=\Gal(F_{\infty}/\Q^{\cyc})$ and $\Ga = \Gal(\Q^{\cyc}/\Q)$. Note
that $G$ and $H$ are compact $p$-adic group with no $p$-torsion but
they are clearly not pro-$p$. For $n\geq 1$, let $\rho_n$ denote the
representation of $G$ obtained by inducing any character of exact
order $p^n$ of $\Gal(\Q(\mu_{p^n}, m^{p^{-n}})/\Q(\mu_{p^n}))$ to
$\Gal(\Q(\mu_{p^n}, m^{p^{-n}})/\Q)$. This representation is defined
over $\Q$ and is irreducible of dimension $p^{n-1}(p-1)$. We may
view $\rho_n$ as an Artin representation $G\lra
\mathrm{GL}_{p^{n-1}(p-1)}(\Op)$. Now if $M$ is a
$\Op\ps{G}$-module, we define $\mathrm{tw}_{\rho_n}(M)$ to be the
$\Op$-module $W_{\rho_n}\ot_{\Zp}M$ with $G$ acting diagonally. Here
$W_{\rho_n}$ is a free $\Op$-module of rank $p^{n-1}(p-1)$ which
realizes $\rho_n$. We are in the position to state the following
theorem.

\bt \label{Akashi twist}
 Retain the assumptions of Theorem \ref{Akashi congruent} and notation as above.
 Furthermore, assume that $X(A/F_{\infty})$ and $X(B/F_{\infty})$
 lie in $\M_H(G)$ $($see Section \ref{char section} for definition$)$.
Then $Ak_{H} \big(\mathrm{tw}_{\rho_n}X(A/F_{\infty})\big)$ is a
unit in $\Op\ps{\Ga}$ if and only if $Ak_{H}
\big(\mathrm{tw}_{\rho_n}X(B/F_{\infty})\big)$ is a unit in
$\Op\ps{\Ga}$. \et

\bpf
 Let $G_n$ and $G_n'$ be the open subgroup of $G$ corresponding to
$\Q(\mu_{p^n}, m^{p^{-n}})$ and $\Q(\mu_{p^n}, m^{p^{-(n-1)}})$
respectively. Write $H_n =H\cap G_n$, $H_n' =H\cap G_n'$ and $\Ga_n
= G_n/H_n = G_n'/H_n'$. Note that $G_n, G_n', H_n$ and $H_n'$ are
pro-$p$ groups without $p$-torsion. By \cite[Theorem A44]{DD}, we
have
\[ Ak_{H_n}X(A/F_{\infty}) = Ak_{H_n'}X(A/F_{\infty})N_{\Ga/\Ga_n}\Big(Ak_{H}\big(\mathrm{tw}_{\rho_n}X(A/F_{\infty})\big)^{p-1}\Big) \]
 in $Q_{\Op}(\Ga_n)/\Op\ps{\Ga_n}^{\times}$, where $Q_{\Op}(\Ga_n)$
is the field of quotient of $\Op\ps{\Ga_n}$ and $N_{\Ga/\Ga_n}$ is
the norm map from $\Op\ps{\Ga}$ to $\Op\ps{\Ga_n}$.

Now if $Ak_{H} \big(\mathrm{tw}_{\rho_n}X(A/F_{\infty})\big)$ is a
unit in $\Op\ps{\Ga}$, then
$N_{\Ga/\Ga_n}\Big(Ak_{H}\big(\mathrm{tw}_{\rho_n}X(A/F_{\infty})\big)^{p-1}\Big)$
is a unit in $\Op\ps{\Ga_n}$. Therefore, we have
\[ Ak_{H_n}X(A/F_{\infty}) = Ak_{H_n'}X(A/F_{\infty}) ~\mbox{mod}~
\Op\ps{\Ga_n}^{\times}. \]
 Write $F_n = \Q(\mu_{p^n}, m^{p^{-n}})$. From now on, identify
$\Op\ps{\Ga_n}\cong \Op\ps{T}$ under a choice of a generator of
$\Ga_n$. Recall that a polynomial $T^n + c_{n-1}T^{n-1} + \cdots
+c_0$ in $\Op[T]$ is said to be a Weierstrass polynomial if $\pi$
divides $c_i$ for every $0\leq i \leq n-1$. Since $X(A/F_{\infty})$
lies in $\M_H(G)$, by restricting scalars, we also have that
$X(A/F_{\infty})$ lies in $\M_{H_n}(G_n)$. Since $G_n$ is pro-$p$,
we may apply \cite[Theorem 3.1]{LimMu} to conclude that
\[
\mu_{\Op\ps{\Ga_n}}\big(X(A/F_n^{\cyc})\big) =
\mu_{\Op\ps{G_n}}\big(X(A/F_{\infty})\big).\]
  By virtue of
Proposition \ref{Selmer fg} and the above equality, we have
 \[Ak_{H_n}X(A/F_{\infty}) =
 \pi^{\mu_{\Op\ps{G_n}}\big(X(A/F_{\infty})\big)}\frac{f}{g}\]
 for some Weierstrass polynomials $f$ and $g$.
 Similarly, we have \[Ak_{H_n'}X(A/F_{\infty}) =
 \pi^{\mu_{\Op\ps{G_n'}}\big(X(A/F_{\infty})\big)}\frac{f'}{g'}\]
 for some Weierstrass polynomials $f'$ and $g'$. Therefore, it
 follows from the above Akashi series relation that
$\pi^{\mu_{\Op\ps{G_n}}\big(X(A/F_{\infty})\big)}fg'$ and
$\pi^{\mu_{\Op\ps{G_n'}}\big(X(A/F_{\infty})\big)}f'g$ generates the
same ideal in $\Op\ps{T}$. By an $\Op$-analogue of \cite[Lemma
5.7]{LimMHG}, this in turn implies that
\[ \mu_{\Op\ps{G_n}}\big(X(A/F_{\infty})\big) =
\mu_{\Op\ps{G_n'}}\big(X(A/F_{\infty})\big)\] and
 \[ \deg f - \deg g = \deg f' - \deg g'.\]
 On the other hand, since $G'$ is a subgroup of $G$ with index $p$, we have
 \[ \mu_{\Op\ps{G_n'}}\big(X(A/F_{\infty})\big) =
p\mu_{\Op\ps{G_n}}\big(X(A/F_{\infty})\big)\] and combining this
with the above, we obtain
\[ \mu_{\Op\ps{G_n}}\big(X(A/F_{\infty})\big) =
\mu_{\Op\ps{G_n'}}\big(X(A/F_{\infty})\big) = 0.\]
   Since
$X(A/F_{\infty})$ is assumed to lie in $\M_{H_n}(G_n)$, it follows
from \cite[Proposition 5.7]{LimMu} that $X(A/F_{\infty})$ is
finitely generated over $\Op\ps{H_n}$. In particular, by the
argument in the proof of \cite[Lemma 5.4]{LimMHG}, we have
\[ \rank_{\Op\ps{H_n}}\big(X(A/F_{\infty})\big) = \deg f- \deg g \]
Similarly, we have
 \[ \rank_{\Op\ps{H'_n}}\big(X(A/F_{\infty})\big) = \deg f'- \deg g'. \]
  Since the right hand side of the two equalities are the same, we have
 \[ \rank_{\Op\ps{H_n}}\big(X(A/F_{\infty})\big) = \rank_{\Op\ps{H'_n}}\big(X(A/F_{\infty})\big)
 .\] On the other hand, since $H_n'$ is a subgroup of $H$ of index $p$, we
 also have \[
 p\rank_{\Op\ps{H_n}}\big(X(A/F_{\infty})\big) =
 \rank_{\Op\ps{H'_n}}\big(X(A/F_{\infty})\big).\]
  Hence we conclude that
\[ \rank_{\Op\ps{H_n}}\big(X(A/F_{\infty})\big) = 0.\] Therefore, we have
shown that $X(A/F_{\infty})$ is a finitely generated torsion
$\Op\ps{H_n}$-module. By Proposition \ref{Akashi congruent lemma},
this in turn implies that $X(B/F_{\infty})$ is also a finitely
generated torsion $\Op\ps{H_n}$-module. Since $F_{\infty}$ is
clearly not totally real and $G_n$ is pro-$p$, we may apply
Proposition \ref{algebra nonzero} to conclude that
$X(B/F_{\infty})=0$. It follows from this observation that
$\mathrm{tw}_{\rho_n}X(B/F_{\infty})=0$ which in turn implies that
$Ak_{H} \big(\mathrm{tw}_{\rho_n}X(B/F_{\infty})\big)$ is a unit in
$\Op\ps{\Ga}$. This completes the proof of the theorem. \epf

The astute reader will notice that the proof given above can also
establish the following proposition (compare with \cite[Theorem
A.13]{DD}). Again, we note that we do not assume that
$X(A/F_{\infty})$ has no non-zero pseudo-null $\Op\ps{G}$-submodule
(but we require $F_{\infty}$ to be not totally real).

\bp \label{Akashi twist prop}
 Retain the setting of Theorem \ref{Akashi twist}.
Then $X(A/F_{\infty})=0$  if and only if there exists $n$ such that
$Ak_{H} \big(\mathrm{tw}_{\rho_n}X(A/F_{\infty})\big)$ is a unit in
$\Op\ps{\Ga}$. \ep

\section{Comparing characteristic elements of Selmer groups}
\label{char section}

In this section, we will prove our main results concerning with the
characteristic elements of the Selmer groups. We need to introduce
some further notion and notation. Let
\[ \Si = \{\,s\in \Op\ps{G}~ | ~\Op\ps{G}/\Op\ps{G}s
~\mbox{is a finitely generated $\Op\ps{H}$-module} \}.\]
 By
\cite[Theorem 2.4]{CFKSV}, $\Si$ is a left and right Ore set
consisting of non-zero divisors in $\Op\ps{G}$. Set $\Si^* =
\cup_{n\geq 0}\pi^n\Si$. It follows from \cite[Proposition
2.3]{CFKSV} that a finitely generated $\Op\ps{G}$-module $M$ is
annihilated by $\Si^*$ if and only if $M/M(\pi)$ is finitely
generated over $\Op\ps{H}$. We will denote $\M_H(G)$ to be the
category of all finitely generated $\Op\ps{G}$-modules which are
$\Si^*$-torsion. For data arising from abelian varieties and modular
forms, it has been conjectured that the dual Selmer groups
associated to these data lie in this category (see \cite{CFKSV,
CS12, Su}).

By the discussion in \cite[Section 3]{CFKSV}, we have the following
exact sequence
\[ K_1(\Op\ps{G}) \lra K_1(\Op\ps{G}_{\Si^*})
\stackrel{\partial_G}{\lra} K_0(\M_H(G))\lra 0 \] of $K$-groups. For
each $M$ in $\M_H(G)$, we define a \textit{characteristic element}
for $M$ to be any element $\xi_M$ in $K_1(\Zp\ps{G}_{\Si^*})$ with
the property that
\[\partial_G(\xi_M) = -[M].\]

Let $\rho:G\lra GL_m(\Op_{\rho})$ denote a continuous group
representation (not necessarily an Artin representation), where
$\Op'= \Op_{\rho}$ is the ring of integers of some finite extension
of $K$. For $g\in G$, we write $\bar{g}$ for its image in $\Ga=G/H$.
We define a continuous group homomorphism
 \[G \lra M_d(\Op')\ot_{\Op}\Op\ps{\Ga}, \quad g\mapsto \rho(g)\ot
 \bar{g}. \] By \cite[Lemma 3.3]{CFKSV}, this in turn induces a map
 \[ \Phi_{\rho}: K_1(\Op\ps{G}_{\Si^*})\lra Q_{\Op'}(\Ga)^{\times}, \]
where $Q_{\Op'}(\Ga)$ is the field of fraction of $\Op'\ps{\Ga}$.
Let $\varphi: \Op'\ps{\Ga}\lra \Op'$ be the augmentation map and
denote its kernel by $\mathfrak{p}$. One can extend $\varphi$ to a
map $\varphi : \Op'\ps{\Ga}_{\mathfrak{p}}\lra K'$, where $K'$ is
the field of fraction of $\Op'$. Let $\xi$ be an arbitrary element
in $K_1(\Op\ps{G}_{\Si^*})$. If $\Phi_{\rho}(\xi)\in
\Op'\ps{\Ga}_{\mathfrak{p}}$, we define $\xi(\rho)$ to be
$\varphi(\Phi_{\rho}(\xi))$. If $\Phi_{\rho}(\xi)\notin
\Op'\ps{\Ga}_{\mathfrak{p}}$, we set $\xi(\rho)$ to be $\infty$.

We will write $\al$ for the natural map
\[\Op\ps{G}_{\Si^*}^{\times}\lra K_1(\Op\ps{G}_{\Si^*}).\]
We can now state the following result which will prove
\cite[Conjecture 4.8 Case 4]{CFKSV} for the characteristic elements
attached to our Selmer groups. We mention that this result refines
the previous result of the author \cite[Proposition 6.3]{LimMHG},
where he proved the same result but under the extra assumption that
$G=\Gal(F_{\infty}/F)$ is pro-$p$.

\bp \label{conjecture 4.8} Let $F_{\infty}$ be a strongly admissible
$p$-adic Lie extension of $F$, and assume that $F_{\infty}$ is not
totally real. Suppose that $\mathbf{(Fin_{\infty})}$ and
$\mathbf{(Cong)}$ are satisfied. Also, suppose that
$X(A/F_{\infty})\in \M_H(G)$. Let $\xi_A$ be a characteristic
element of $X(A/F_{\infty})$. Then the following statements are
equivalent.

\begin{enumerate}
\item[$(a)$] $\xi_A \in \al(\Op\ps{G}^{\times})$, where $\al$ is
the map $\Op\ps{G}_{\Si^*}^{\times}\lra K_1(\Op\ps{G}_{\Si^*})$.

\item[$(b)$] $\xi_A(\rho)$ is finite and lies in $\Op_{\rho}^{\times}$ for every continuous
group representation $\rho$ of $G$.

\item[$(c)$] $\Phi_{\rho}(\xi_A) \in \Op_{\rho}\ps{\Ga}^{\times}$ for every continuous
group representation $\rho$ of $G$.

\item[$(d)$] $\Phi_{\rho}(\xi_A) \in \Op_{\rho}\ps{\Ga}^{\times}$ for every Artin
 representation $\rho$ of $G$.

\item[$(e)$] There exists an open normal pro-$p$ subgroup $G'$ of
$G$ such $Ak_{H'}\big(X(A/F_{\infty})\big) \in
\Op\ps{\Ga'}^{\times}$. Here $H' = H\cap G'$ and $\Ga' = G'/H'$.
\end{enumerate}
\ep

\bpf
 By \cite[Lemma 4.9]{CFKSV}, we have the implications
(a)$\Rightarrow$(b)$\Leftrightarrow$(c)$\Rightarrow$(d). It remains
to show (d)$\Rightarrow$(e)$\Rightarrow$(a). We first establish the
implication (d)$\Rightarrow$(e).
 Let $G'$ be an open normal subgroup of $G$ with the property that
 $G'$ is pro-$p$. Write $H'= H\cap G'$ and $\Ga' = G'/H'$. Let
 $\Op'$ denote the ring of integers of a finite extension $K'$ of
 $\Qp$ such that $K'$ contains $K$ and all absolutely irreducible
 representations of $\Delta= G/G'$ can be realized over $K'$.
 Denote $\hat{\Delta}$ to be the set of all irreducible representations
 of $\Delta$. Write $M_{\Op'} = M\ot_{\Op}\Op'$. If $M$ is a
$\Op\ps{G}$-module and $\rho\in\hat{\Delta}$ with dimension
$n_{\rho}$, we define $\mathrm{tw}_{\rho}(M)$ to be the
$\Op'$-module $W_{\rho}\ot_{\Op'}M_{\Op'}$ with $G$ acting
diagonally. Here $W_{\rho}$ is a free $\Op'$-module of rank
$n_{\rho}$ realizing $\rho$.
 Then by
 \cite[Theorem A44]{DD}, we have
 \[ Ak_{H'} \big(X(A/F_{\infty})_{\Op'}\big)^{[\Ga : \Ga']}
  = N_{\Ga/\Ga'}\Big( \prod_{\rho\in\widehat{\Delta}}
  Ak_{H} \big(\mathrm{tw}_{\rho}X(A/F_{\infty})\big)^{n_{\rho}} \Big)~\mathrm{mod}\,\Op'\ps{\Ga}^{\times}, \]
  where $N_{\Ga/\Ga'}$ is the norm map from $\Op\ps{\Ga}$ to
  $\Op\ps{\Ga'}$ and $n_{\rho}$ is the dimension of $\rho$.
On the other hand, it follows from \cite[Lemma 3.7]{CFKSV} that one
has
\[ \Phi_{\rho}(\xi_A) =
Ak_{H} \big(\mathrm{tw}_{\hat{\rho}}X(A/F_{\infty})\big)
~\mathrm{mod}\,\Op'\ps{\Ga}^{\times},\] where $\hat{\rho}$ is the
contragredient of $\rho$. Now if statement (d) holds, it then
follows from combining the above two observations that
\[ Ak_{H'}\big(X(A/F_{\infty})_{\Op'}\big) \in \Op'\ps{\Ga}^{\times}.\]
It is now an easy exercise to deduce from this that
\[Ak_{H'}\big(X(A/F_{\infty})\big) \in
\Op\ps{\Ga}^{\times}.\] This proves the implication
(d)$\Rightarrow$(e)

We now prove (e)$\Rightarrow$(a). Suppose that statement (e) holds.
Let $F'$ be the fixed field of $G'$ as given by (e). Since we are
assuming that $X(A/F_{\infty})\in \M_H(G)$, we may apply
\cite[Proposition 2.5]{CS12} to conclude that $X(A/L^{\cyc})$ is
$\Op\ps{\Ga_L}$-torsion for every finite extension $L$ of $F$
contained in $F_{\infty}$, where $\Ga_L = \Gal(L^{\cyc}/L)$.
Combining this with our $\mathbf{(Fin_{\infty})}$ assumption, it
follows from a similar argument to that in \cite[Proposition 3.3(i),
Corollary 3.4]{LimMHG} that $H^2(G_S(F_{\infty}),A)=0$ and
$\la_{A/F_{\infty}}$ is surjective. Therefore, we may apply
Propositions \ref{algebra prop} and \ref{Selmer fg} to conclude that
$X(A/F_{\infty})$ is a finitely generated torsion
$\Op\ps{H'}$-module. Since $F_{\infty}$ is not totally real and $G'$
is pro-$p$, it follows from Lemma \ref{algebra nonzero}(a) that
$X(A/F_{\infty}) =0$. In particular, this implies that
$\partial_G(\xi_A)=0$. Consequently, it follows from the exact
sequence of $K$-groups
\[ K_1(\Op\ps{G}) \lra K_1(\Op\ps{G}_{\Si^*})
\stackrel{\partial_G}{\lra} K_0(\M_H(G))\lra 0 \] that there exists
an element in $K_1(\Op\ps{G})$ which maps to $\xi_A$. On the other
hand, since $\Op\ps{G}$ is semi-local, we have that
$\Op\ps{G}^{\times}$ maps onto $K_1(\Op\ps{G})$, thus yielding (a).
 \epf

We mention that one can also establish the following proposition
which will refine \cite[Proposition 6.2]{LimMHG}. As the proof is
very similar to the preceding proposition, we will omit it.

 \bp \label{conjecture 4.8'} Let $M\in \M_H(G)$. Suppose that
 $M$ contain no nontrivial pseudo-null $\Op\ps{G}$-submodules.
 Let $\xi_M$ be a characteristic
element of $M$. Then the following statements are equivalent.

\begin{enumerate}
\item[$(a)$] $\xi_M \in \al(\Op\ps{G}^{\times})$, where $\al$ is
the map $\Op\ps{G}_{\Si^*}^{\times}\lra K_1(\Op\ps{G}_{\Si^*})$.

\item[$(b)$] $\xi_M(\rho)$ is finite and lies in $\Op_{\rho}^{\times}$ for every continuous
group representation $\rho$ of $G$.

\item[$(c)$] $\Phi_{\rho}(\xi_M) \in \Op_{\rho}\ps{\Ga}^{\times}$ for every continuous
group representation $\rho$ of $G$.

\item[$(d)$] $\Phi_{\rho}(\xi_M) \in \Op_{\rho}\ps{\Ga}^{\times}$ for every Artin
 representation $\rho$ of $G$.
 \end{enumerate}
\ep

 We can now prove our final result which compares the
 characteristic elements of the Selmer groups of two congruent data.

\bt \label{char congruent}
 Let $F_{\infty}$ be a strongly admissible $p$-adic Lie extension of $F$.
 Assume that $F$ is not totally real.
 Suppose that all the hypothesis in Proposition \ref{Selmer fg} are valid
 for $A$ and $B$, and that
 $\mathbf{(Cong)}$ is satisfied.
Furthermore, suppose that for every $v\notin
 S_2$, we have
 \[\dim_k \big(C_v(F^{\cyc}_w)/\pi\big) = \dim_k\big(D_v(F^{\cyc}_w)/\pi\big)\]
  for every prime $w$ of $F^{\cyc}$ above $v$.

Then $\xi_A\in \al(\Op\ps{G}^{\times})$ if and only if $\xi_{B}\in
\al(\Op\ps{G}^{\times})$. \et

\bpf
 Suppose that $\xi_A\in \al(\Op\ps{G}^{\times})$. Then by
Proposition \ref{conjecture 4.8}, there exists an open normal
pro-$p$ subgroup $G'$ of $G$ such $Ak_{H'}\big(X(A/F_{\infty})\big)
\in \Op\ps{\Ga}^{\times}$, where $H' = H\cap G'$. By Theorem
\ref{Akashi congruent}, this in turn implies that
$Ak_{H'}\big(X(B/F_{\infty})\big) \in \Op\ps{\Ga}^{\times}$. We may
now apply Proposition \ref{conjecture 4.8} again to conclude that
$\xi_{B}\in \al(\Op\ps{G}^{\times})$.
 \epf

\begin{ack}
   We like to thank the anonymous referee for pointing out
    some mistakes
and valuable comments which improved the clarity of the paper.
            \end{ack}


\footnotesize
\begin{thebibliography}{9999999}

\bibitem[Ch]{Ch} A. Chandrakant Sharma, Non existence of finite $\La$-submodules
of dual Selmer groups over a cyclotomic $\Zp$-extension, \textit{J.
Ramanujan Math. Soc.} \textbf{24} (2009), No. 1, 75--85.

\bibitem[Ch2]{Ch09} A. Chandrakant Sharma, Iwasawa invariants for the False-Tate
extension and congruences between modular forms, \emph{J. Number
Theory} 129 (2009) 1893--1911.

\bibitem[CFKSV]{CFKSV} J.\ Coates, T.\ Fukaya, K.\ Kato, R.\ Sujatha
and O.\ Venjakob, The $GL_2$ main conjecture for elliptic curves
without complex multiplication, \textit{Publ. Math. IHES}
\textbf{101} (2005), 163--208.

\bibitem[CG]{CG}  J. Coates and R. Greenberg, Kummer theory for abelian varieties over
local fields, \emph{Invent. Math.} 124 (1996) 129--174.

\bibitem[CSS]{CSS} J. Coates, P. Schneider and R. Sujatha, Links between cyclotomic
and $GL_2$ Iwasawa theory. Kazuya Kato's fiftieth birthday.
\textit{Doc. Math.} 2003, Extra Vol., 187--215 (electronic).

\bibitem[CS]{CS12} J.\ Coates and R.\ Sujatha, On the
$\mathfrak{M}_H(G)$-conjecture, in \textit{Non-abelian fundamental
groups and Iwasawa theory}, 132--161, London Math. Soc. Lecture Note
Ser., \textbf{393}, Cambridge Univ. Press, Cambridge, 2012.

\bibitem[DD]{DD} T. Dokchitser and V. Dokchitser,
Computations in non-commutative Iwasawa theory. With an appendix by
J. Coates and R. Sujatha. \textit{Proc. Lond. Math. Soc. (3)}
\textbf{94} (2007), no. 1, 211--272.

\bibitem[EPW]{EPW} M. Emerton, R. Pollack and T. Weston,
Variation of Iwasawa invariants in Hida families, \emph{Invent.
Math.} 163 (2006) 523--580.

\bibitem[GW]{GW} K. R. Goodearl and R. B. Warfield, \textit{An
introduction to non-commutative Noetherian rings}, London Math. Soc.
Stud. Texts \textbf{61}, Cambridge University Press, 2004.

\bibitem[G89]{G89} R. Greenberg, Iwasawa theory for $p$-adic representations, in
\emph{Algebraic Number Theory--in honor of K. Iwasawa}, ed. J.
Coates, R. Greenberg, B. Mazur and I. Satake, Adv. Std. in Pure
Math. 17, 1989, pp. 97--137.

\bibitem[G94]{Gr94} R. Greenberg, Iwasawa theory for $p$-adic deformations
of motives, \emph{Proc. Sympos. Pure Math.} 55 (Part 2) (1994)
193--223.

\bibitem[GV]{GV} R. Greenberg and V. Vatsal, On the Iwasawa invariants of elliptic
curves, \emph{Invent. Math.} \textbf{142} (2000) 17--63.

\bibitem[Ha]{Ha}  Y. Hachimori, Iwasawa $\lambda$-invariants and
congruence of Galois representations, \emph{J. Ramanujan Math. Soc.}
\textbf{26}(2) (2011) 203--217.

\bibitem[HM]{HM} Y. Hachimori and K. Matsuno, An analogue of Kida's
formula for the Selmer groups of elliptic curves, \textit{J.
Algebraic Geom.} \textbf{8} (1999), no. 3, 581--601.

\bibitem[HO]{HO} Y.\ Hachimori and T.\ Ochiai, Notes on
non-commutative Iwasawa theory, \textit{Asian J. Math.} \textbf{14}
(2010), no. 1, 11--17.

\bibitem[Ho]{Ho} S. Howson, Euler characteristic as invariants of
Iwasawa modules, \textit{Proc. London Math. Soc. (3)} \textbf{85}
(2002), no. 3, 634--658.

\bibitem[Im]{Im} H. Imai, A remark on the rational points of abelian
varieties with values in cyclotomic $\Zp$-extensions, \emph{Proc.
Japan Acad.} 51 (1975) 12--16.

\bibitem[Lam]{Lam} T. Y. Lam, \textit{Lectures on Modules and Rings},
Grad. Texts in Math. \textbf{189},
Springer-Verlag, New York, 1999.

\bibitem[Lim]{LimMHG}  M. F. Lim, A remark on the $\mathfrak{M}_H(G)$-conjecture
and Akashi series, \textit{Int. J. Number Theory} \textbf{11}
(2015), No. 1, 269--297.

\bibitem[LimMu]{LimMu} M. F. Lim, Comparing the Selmer group of a $p$-adic
representation and the Selmer group of the Tate dual of the
representation, arXiv:1405.5289 [math.NT].

\bibitem[Mat]{Mat} K. Matsuno, Finite $\La$-submodules of Selmer groups of abelian
varieties over cyclotomic $\Zp$-extensions, \textit{J. Number
Theory} \textbf{99} (2003), no. 2, 415--443.

\bibitem[MW]{MW} B. Mazur and A. Wiles, On $p$-adic analytic families of Galois
representations, \textit{Compo. Math.} \textbf{59} (1986) 231--264.

\bibitem[Neu]{Neu} A. Neumann, Completed group algebras without zero divisors, \textit{Arch.
Math.} \textbf{51} (1988), no. 6, 496--499.

\bibitem[Oc]{Oc} Y.\ Ochi, A note on Selmer groups of abelian
varieties over the trivializing extensions, \textit{Proc. Amer.
Math. Soc.} \textbf{134} (2006), no. 1, 31--37.

\bibitem[OcV]{OcV02} Y.\ Ochi and O. Venjakob, On the structure of
Selmer groups over $p$-adic Lie extensions, \textit{J. Algebraic
Geom.} \textbf{11} (2002), no. 3, 547--580.

\bibitem[Sh]{Sh} S. Shekhar, Euler characteristics of $\La$-adic
forms over Kummer extensions, \textit{Int. J. Number Theory}
\textbf{10} (2014), No. 2, 401--420.

\bibitem[SS]{SS} S. Shekhar and R. Sujatha, On the structure of
Selmer groups of $\La$-adic deformations over $p$-adic Lie
extensions, \textit{Doc. Math.} \textbf{17} (2012), 573--606.

\bibitem[SS2]{SS14} S. Shekhar and R. Sujatha, Euler characteristic and
congruences of elliptic curves, \textit{M\"unster J. of Math.} 7
(2014), 327--343.

\bibitem[Su]{Su} R. Sujatha, Iwasawa theory and modular forms,
\textit{Pure Appl. Math. Q.} \textbf{2} (2006), no. 2, 519--538.

\bibitem[V]{V02} O. Venjakob, On the structure theory of the Iwasawa algebra
of a $p$-adic Lie group, \textit{J. Eur. Math. Soc.} \textbf{4}
(2002), no. 3, 271--311.

\end{thebibliography}
\end{document}